\documentclass[twoside,12pt]{article}
\usepackage{latexsym}
\usepackage{amsmath}
\usepackage{amssymb}
\usepackage{cases}
\usepackage{ascmac}
\usepackage{mathrsfs}
\usepackage{amsthm}
\usepackage{cite}
\usepackage[usenames,dvipsnames]{color}

\title{Abstract approach of degenerate parabolic equations with dynamic boundary conditions}

\author{Takeshi Fukao\\
Department of Mathematics, Faculty of Education\\
Kyoto University of Education\\
1~Fujinomori, Fukakusa, Fushimi-ku, Kyoto~612-8522 Japan\\
E-mail: \texttt{fukao@kyokyo-u.ac.jp}\\
\and \\ Taishi Motoda\\
Graduate School of Education, 
Kyoto University of Education\\
1~Fujinomori, Fukakusa, Fushimi-ku, Kyoto~612-8522 Japan\\
E-mail: \texttt{motoda.math@gmail.com}}
 
\date{}

\newcommand\testopari{\sc Takeshi Fukao and Taishi Motoda}
\newcommand\testodispari{\sc Degenerate parabolic equations with dynamic boundary conditions}
\begin{document}

\maketitle

\begin{abstract}

An initial boundary value problem of the nonlinear diffusion equation with 
a dynamic boundary condition is treated. 
The existence problem of the initial-boundary value problem is discussed. 
The main idea of the proof is an 
abstract approach from the evolution equation governed by the subdifferential.
To apply this, the setting of suitable function spaces, more precisely 
the mean-zero function spaces, is important. 
In the case of a dynamic boundary condition, the total mass, which is the sum of 
volumes in the bulk and on the boundary, is a point of emphasis. 
The existence of a weak solution is proved on this basis.

\vspace{2mm}

\vspace{2mm}

\end{abstract}

\section{Introduction}
\setcounter{equation}{0}

We consider the initial boundary value problem of the nonlinear diffusion equation 
{\rm (P)}, comprising 
\begin{gather}
	\displaystyle \frac{\partial u}{\partial t}-\Delta \xi =f, \quad \xi \in \beta (u) 
	\quad {\rm in}~Q:=(0, T)\times \Omega, 
	\label{equ1}
	\\
	\displaystyle \xi _{|_{\Gamma }}=\xi _{\Gamma }, 
	\quad \frac{\partial u_{\Gamma }}{\partial t}+
	\partial _{\boldsymbol{\nu }}\xi -\Delta _{\Gamma }\xi _{\Gamma }
	=f_{\Gamma },
	\quad \xi _{\Gamma }\in \beta (u_{\Gamma }) 
	\quad {\rm on}~\Sigma :=(0, T)\times \Gamma, 
	\label{equ2} 
	\\
	u(0)=u_{0} 
	\quad {\rm in}~\Omega, 
	\quad u_{\Gamma }(0)=u_{0\Gamma }
	\quad {\rm on}~\Gamma, 
	\label{equ3}
\end{gather}
where $0< T< +\infty$, $\Omega $ is a bounded domain of $\mathbb{R}^{d}$ $(d=2,3)$ 
with smooth boundary $\Gamma :=\partial \Omega $, 
$\xi _{|_{\Gamma }}$ stands for the trace of $\xi $ to $\Gamma $, 
$\partial _{\boldsymbol{\nu}}$ is the outward normal derivative on $\Gamma $, 
$\Delta $ is the {L}aplacian, $\Delta _{\Gamma }$ is the {L}aplace--{B}eltrami operator 
(see, e.g., \cite{Gri09}), and
$f: Q\to \mathbb{R}$, 
$f_{\Gamma }: \Sigma \to \mathbb{R}$, 
$u_{0}: \Omega \to \mathbb{R}$, and 
$u_{0\Gamma }: \Gamma \to \mathbb{R}$ are given data. 
Moreover,
$\beta : \mathbb{R}\to 2^{\mathbb{R}}$, a maximal monotone graph, 
characterizes the first and second equations of {\rm (P)}
as the degenerate parabolic system. Indeed, 
by choosing various types of $\beta $ given later, 
{\rm (P)} will be various types of the degenerate parabolic system; e.g., 
{\rm (P)} can be the {S}tefan problem, porous media equation, 
or fast diffusion equation (see, e.g., \cite{CF16}). 
In particular, we allow $\beta $ to be multivalued 
because we are also interested in the Hele--Shaw profile; 
more precisely, 
$\beta:=\partial I_{[0, 1]}$, 
the subdifferential of the indicator function $I_{[0,1]}$ on interval $[0,1]$. 
In this paper, we treat a modified version of 
the Hele--Shaw profile.

In terms of the well-posedness of {\rm (P)}, an early result for the 
{S}tefan problem was given \cite{Dam77}. 
For this result, an abstract theory of the evolution equation in {H}ilbert space was 
applied. 
There are also treatments of {\rm (P)} \cite{Bar10}; more precisely, there 
are two major approaches named the {H}ilbert space approach and $L^1$ approach. 
Results obtained using the {H}ilbert space approach have been presented \cite{Aik93, Aik95, Aik96} 
related to the {S}tefan problem with a dynamic boundary condition, 
\cite{Ken90, DK99} for a wider degenerate parabolic equation. 
Results obtained using the $L^1$ approach have been reported \cite{IK02, AMTI06, Igb07}. 
We refer to \eqref{equ2}, which includes a time derivative, as the dynamic boundary condition. 
Asymptotic analysis 
of the {C}ahn--{H}iliard equation has recently been performed \cite{Fuk16, Fuk16b}. 
In this treatment, if we choose different values of $\beta $ between \eqref{equ1} and \eqref{equ2}, namely 
$\beta $ and $\beta _\Gamma $, then we need a domination assumption \cite[p.419, {\rm (A6)}]{CF15}. 
We improve this assumption in Section~4 of the present paper. 
In the cited studies, 
the important point is the setting of function spaces, where 
the total mass is zero. 
This property arises from the dynamic boundary condition 
(see also \cite{Ken90, CF16, KY17} for the setting of the {N}eumann boundary condition).  
In the present paper, to apply the pioneering idea of \cite{Dam77}, 
we use the same setting \cite{Fuk16, Fuk16b} to construct the duality mapping that plays the role of 
diffusion. 
One of the greatest difficulties of the problem is similar to the case of the {N}eumann boundary condition (see, e.g., 
\cite{Ken90, KL05}).

The present paper proceeds as follows. 
Section~2 states the main theorem. 
We first prepare the notation used in this paper and set the 
suitable duality mapping and function spaces. 
We then introduce the definition of the weak solution of {\rm (P)}, 
and give the main theorem.

In Section~3, to apply the 
abstract theory of the evolution equation governed by 
the subdifferential \cite{Bre73}, 
we define the proper lower semicontinuous and convex functional. 
We consider the approximate problem using {M}oreau--{Y}osida regularization. 
We also give characterization lemma for the subdifferential. 
We then deduce uniform estimates of the approximate solutions.
We finally prove the existence of weak solutions by passing to the limit.

In Section~4, we discuss improvements to the assumptions.

A detailed index of sections and subsections follows.

\begin{itemize}
 \item[1.] Introduction
 \item[2.] Main results
\begin{itemize}
 \item[2.1.] Notation
 \item[2.2.] Definition of the solution and main theorem
\end{itemize}
 \item[3.] Approximate problem and uniform estimates
\begin{itemize}
 \item[3.1.] Abstract formulation
 \item[3.2.] Approximate problem for {\rm (P)}
 \item[3.3.] Uniform estimates
 \item[3.4.] Passage to the limit as $\lambda \to 0$
\end{itemize}
 \item[4.] Improvement
 \begin{itemize}
 \item[4.1.] Improvement of the initial condition to a nonzero mean value
 \item[4.2.] Nonlinear diffusions of different $\beta $ and $\beta _\Gamma $
\end{itemize}

\end{itemize}

\section{Main results}
\setcounter{equation}{0}

\subsection{Notation}
We use the spaces 
$H:=L^{2}(\Omega )$, 
$H_{\Gamma }:=L^{2}(\Gamma )$, 
$V:=H^{1}(\Omega )$, 
$V_{\Gamma }:=H^{1}(\Gamma )$ 
with respective standard norms 
$|\cdot |_{H}$, 
$|\cdot |_{H_{\Gamma }}$, 
$|\cdot |_{V}$, 
$|\cdot |_{V_{\Gamma }}$ 
and inner products 
$(\cdot , \cdot )_{H}$, 
$(\cdot , \cdot )_{H_{\Gamma }}$, 
$(\cdot , \cdot )_{V}$, 
$(\cdot , \cdot )_{V_{\Gamma }}$. 
Moreover, we set $\boldsymbol{H}:=H\times H_{\Gamma }$ and
\begin{equation*}
	\boldsymbol{V}:= 
	\bigl\{ 
	\boldsymbol{z}:=(z, z_{\Gamma }) 
	\in V \times V_{\Gamma } \ : \ z_{|_\Gamma }=z_{\Gamma } \ {\rm a.e.~on}~\Gamma 
	\bigr\}. 
\end{equation*}
$\boldsymbol{H}$ and $\boldsymbol{V}$ are then {H}ilbert spaces with inner products
\begin{gather*}
	(\boldsymbol{u}, \boldsymbol{z})_{\boldsymbol{H}}
	:=(u, z)_{H}+(u_{\Gamma }, z_{\Gamma })_{H_{\Gamma }} 
	\quad {\rm for~all~} 
	\boldsymbol{u}:=(u, u_{\Gamma }), \boldsymbol{z}:=(z, z_{\Gamma })\in \boldsymbol{H}, 
	\\
	(\boldsymbol{u}, \boldsymbol{z})_{\boldsymbol{V}}
	:=(u, z)_{V}+(u_{\Gamma }, z_{\Gamma })_{V_{\Gamma }} 
	\quad \quad {\rm for~all~} \boldsymbol{u}
	:=(u, u_{\Gamma }), \boldsymbol{z}:=(z, z_{\Gamma })\in \boldsymbol{V}.
\end{gather*}
Note that 
$\boldsymbol{z} \in \boldsymbol{V}$ 
implies that the second component $z_{\Gamma }$ of $\boldsymbol{z}$ 
is equal to the trace of the first component $z$ of $\boldsymbol{z}$ on $\Gamma $, and 
$\boldsymbol{z} \in \boldsymbol{H}$ implies that 
$z\in H$ and $z_{\Gamma }\in H_{\Gamma }$ are independent. 
Throughout this paper,
we use the bold letter $\boldsymbol{u}$ to represent 
the pair corresponding to the letter; i.e., $\boldsymbol{u}:=(u, u_{\Gamma })$.

Let $m:\boldsymbol{H}\to \mathbb{R}$ be the special mean function defined by 
\begin{equation*}
	m(\boldsymbol{z})
	:=\frac{1}{|\Omega |+|\Gamma |}
	\left\{\int_{\Omega }zdx+\int_{\Gamma }z_{\Gamma }d\Gamma \right\}
	\quad {\rm for~all~} \boldsymbol{z} \in \boldsymbol{H},
\end{equation*}
where $|\Omega |:=\int_{\Omega }1dx, |\Gamma |:=\int_{\Gamma }1d\Gamma $. 
We then define 
$\boldsymbol{H}_{0}:=\{\boldsymbol{z}\in \boldsymbol{H} : m(\boldsymbol{z})=0\}$, $\boldsymbol{V}_{0}:=\boldsymbol{V}\cap \boldsymbol{H}_{0}$. 
Moreover, 
$\boldsymbol{V}^{*}, \boldsymbol{V}_{0}^{*}$ denote the dual spaces of $\boldsymbol{V}, \boldsymbol{V}_{0}$, respectively; 
the duality pairing between $\boldsymbol{V}_{0}^{*}$ and $\boldsymbol{V}_{0}$ is denoted 
$\langle  \cdot , \cdot \rangle _{\boldsymbol{V}_{0}^{*}, \boldsymbol{V}_{0}}$. 
We define the norm of $\boldsymbol{H}_{0}$
by $|\boldsymbol{z}|_{\boldsymbol{H}_{0}}:=|\boldsymbol{z}|_{\boldsymbol{H}}$ 
for all $\boldsymbol{z}\in \boldsymbol{H}_{0}$. 
We now use the 
bilinear form $a(\cdot , \cdot ):\boldsymbol{V}\times \boldsymbol{V}\to \mathbb{R}$, defined by
\begin{equation*} 
	a(\boldsymbol{u}, \boldsymbol{z})
	:= \int_{\Omega } \nabla u\cdot \nabla zdx
	+ \int_{\Gamma } \nabla _{\Gamma } u_{\Gamma }\cdot \nabla _{\Gamma }z_{\Gamma }d\Gamma 
	\quad {\rm for~all~} \boldsymbol{u}, \boldsymbol{z}\in \boldsymbol{V}. 
\end{equation*}
Then, for all $\boldsymbol{z}\in \boldsymbol{V}_{0}$, 
$|\boldsymbol{z}|_{\boldsymbol{V}_{0}}
:=\sqrt{a(\boldsymbol{z}, \boldsymbol{z})}$ is the norm of $\boldsymbol{V}_0$. 
Also, for all $\boldsymbol{z}\in \boldsymbol{V}_{0}$, 
we let $\boldsymbol{F}:\boldsymbol{V}_{0}\to \boldsymbol{V}_{0}^{*}$ be 
the duality mapping defined by
\begin{equation*}
	\langle \boldsymbol{F}\boldsymbol{z}, 
	\tilde{\boldsymbol{z}}
	\rangle _{\boldsymbol{V}_{0}^{*}, \boldsymbol{V}_{0}}
	:=a(\boldsymbol{z}, \tilde{\boldsymbol{z}}) 
	\quad {\rm for~all~}
	\tilde{\boldsymbol{z}}\in \boldsymbol{V}_{0}.
\end{equation*}
Then from the {P}oincar\'e--{W}irtinger inequality, 
there exists a positive constant $c_{\rm P}$ such that 
\begin{equation}
	\label{PW0}
	 |z|_{\boldsymbol{V}}^{2} 
	\le c_{\rm P} \left\{ 
	a(\boldsymbol{z}, \boldsymbol{z})
	+\bigl| m(\boldsymbol{z}) \bigr|^2
	\right\} 
	\quad {\rm for~all~} 
	\boldsymbol{z}\in \boldsymbol{V}.
\end{equation}
Moreover, we define the inner product of $\boldsymbol{V}_{0}^{*}$ by
\begin{equation*}
	(\boldsymbol{z}^{*}, \tilde{\boldsymbol{z}}^{*})_{\boldsymbol{V}_{0}^{*}}
	:=\langle  \boldsymbol{z}^{*}, \boldsymbol{F}^{-1}\tilde{\boldsymbol{z}}^{*}
	\rangle _{\boldsymbol{V}_{0}^{*}, \boldsymbol{V}_{0}} 
	\quad 
	{\rm for~all~} 
	\boldsymbol{z}^{*}, \tilde{\boldsymbol{z}}^{*}\in \boldsymbol{V}_{0}^{*}.
\end{equation*}
We have
$\boldsymbol{V}_{0}\hookrightarrow \hookrightarrow \boldsymbol{H}_{0}\hookrightarrow \hookrightarrow \boldsymbol{V}_{0}^{*}$, where 
``$\hookrightarrow \hookrightarrow $'' stands for 
compact embedding (see \cite[Lemmas A and B]{CF15}). 
One of the essential ideas of the present paper is the setting of the function space $\boldsymbol{V}_0$ and the duality mapping 
$\boldsymbol{F}$ that plays the role of diffusion, as in \cite{Dam77}.

\subsection{Definition of the solution and main theorem}

In this subsection, we define our solution for the 
initial-boundary value problem \eqref{equ1}--\eqref{equ3}, named by {\rm (P)}, and then 
state the main theorem.

\paragraph{Definition 2.1.}
{\it The quadruplet 
$(u, u_{\Gamma }, \xi , \xi _{\Gamma })$ 
is called the weak solution of {\rm (P)} if 
\begin{gather}
	u\in H^{1}(0, T;V^{*})\cap L^{\infty }(0, T;H), 
	\quad u_{\Gamma }\in H^{1}(0, T;V_{\Gamma }^{*})\cap L^{\infty }(0, T;H_{\Gamma }), 
	\label{uug}\\
	\xi \in L^{2}(0, T;V), 
	\quad \xi_\Gamma  \in L^{2}(0, T;V_{\Gamma }), 
	\label{xixig} \nonumber \\
	\xi \in \beta (u) \quad {\it a.e.~in~} Q, 
	\label{xi} \nonumber \\
	\xi _{\Gamma }\in \beta (u_{\Gamma }), 
	\quad \xi _{|_{\Gamma }}=\xi _{\Gamma } 
	\quad {\it a.e.~on~}\Sigma 
	\label{maximal}
\end{gather}
satisfying
\begin{gather}
	\bigl\langle 
	u'(t), z
	\bigr\rangle _{V^{*}, V}+
	\bigl\langle  u'_{\Gamma }(t), z_{\Gamma }
	\bigr\rangle _{V_{\Gamma }^{*}, V_{\Gamma }}
	+\int_{\Omega }\nabla \xi (t)\cdot \nabla zdx
	+\int_{\Gamma }\nabla _{\Gamma }\xi _{\Gamma }(t)\cdot 
	\nabla _{\Gamma }z_{\Gamma }d\Gamma 
	\nonumber \\
	= 
	\int_{\Omega }
	f(t)zdx+\int_{\Gamma }f_{\Gamma }(t)z_{\Gamma }d\Gamma \quad 
	{\it for~all~} \boldsymbol{z}
	:=(z, z_{\Gamma })\in \boldsymbol{V}, 
	\label{vi}
\end{gather}
for a.a.\ $t\in (0, T)$, with}
\begin{equation}
	u(0)=u_{0} \quad {\it a.e.~in}~\Omega , 
	\quad u_{\Gamma }(0)=u_{0\Gamma } 
	\quad {\it a.e.~on}~\Gamma.
	\nonumber \label{ic}
\end{equation}

We assume the following.
\begin{enumerate}
 \item[(A1)] $\beta :\mathbb{R}\to 2^{\mathbb{R}}$ is a maximal monotone graph, which is
the subdifferential $\beta =\partial _{\mathbb{R}}\widehat{\beta }$ 
of some proper lower semicontinuous convex function 
$\widehat{\beta }: \mathbb{R}\to [0, +\infty ]$ 
satisfying $\widehat{\beta }(0)=0$; 
 \item[(A2)] there exist positive constants $c_{1}$, $c_{2}$ such that 
$\widehat{\beta }(r)\geq c_{1}r^{2}-c_{2}$ for all $r\in \mathbb{R}$; 
 \item[(A3)] $\boldsymbol{f}:=(f, f_{\Gamma }) \in L^{2}(0, T; \boldsymbol{H}_0)$; 
 \item[(A4)] $\boldsymbol{u}_{0}:=(u_{0}, u_{0\Gamma })\in \boldsymbol{H}_0, \widehat{\beta }(u_{0})\in L^{1}(\Omega )$ and $\widehat{\beta }(u_{0\Gamma })\in L^{1}(\Gamma )$. 
\end{enumerate}
In particular, (A1) yields $0\in \beta (0)$. 
These assumptions (A1)--(A4) are standard comparing with the literature \cite{CF16, Dam77, DK99, Fuk16, Fuk16b, KL05}. 
Additionally, in the present paper, $\beta $ is modified to a singleton and is similar to a segment far 
from the origin in the following sense.
\begin{enumerate}
 \item[(A5)] There exist constants $c_0$, $M_0>0$, and $c_0'\ge 0$
 such that
\begin{equation}
	\beta (r)
	=
	\begin{cases}
	c_0 r+c_0' & {\rm if~} r \ge  M_0, \\
	c_0 r-c_0' & {\rm if~} r \le  -M_0,
	\end{cases}
	\label{line}
\end{equation}
this implies $D(\beta )=\mathbb{R}$. 
\end{enumerate}

\paragraph{Remark 2.1.} Condition (A5) is a technical yet essential assumption. 
If we can expect that 
the components of $u$ and $u_\Gamma $ of the solution are bounded below by $-M_0$ and above by $M_0$, then 
this modification \eqref{line} is negligible because 
$\beta $ no longer takes these values. 
In the case that we want to treat a maximal monotone graph whose domain is a proper subset of $\mathbb{R}$ (e.g., 
$\partial I_{[0,1]}$), an example of modification is
\begin{equation*}
	\beta (r)=
	\begin{cases}
	r-c_0' & {\rm if~} r < 0, \\
	[-c_0',0] & {\rm if~} r=0, \\
	0 & {\rm if~} 0 < r < 1, \\
	[0,c_0'+1] & {\rm if~} r=1, \\
	r+c_0' & {\rm if~} r > 1.
	\end{cases}
	\quad
	\widehat{\beta }(r)
	=
	\begin{cases}
	\displaystyle \frac{1}{2}r^2-c_0' r & {\rm if~} r < 0, \\
	0 & {\rm if~} 0\le r\le 1, \\
	\displaystyle \frac{1}{2}(r-1)^2+(c_0' +1)(r-1) & {\rm if~} r > 1. \\
	\end{cases}
\end{equation*}
This assumption is used to obtain the uniform boundedness of the total mass.
(cf.\ \cite{Ken90, KL05}). \\

We now give our main theorem.

\paragraph{Theorem 2.1.}
{\it Under assumptions {\rm (A1)}--{\rm (A5)}, 
there exists a unique weak solution to the problem {\rm (P)}. }\\

The continuous dependence of the problem {\rm (P)} is completely the same as that in 
\cite[Theorem~2.2]{Fuk16}, and we therefore omit the proof of the uniqueness in this paper.

\section{Approximate problem and uniform estimates}
\setcounter{equation}{0}

\subsection{Abstract formulation}

We apply the abstract theory of the evolution equation \cite{Bre73} to prove the main theorem, following on from the essential idea of \cite{Dam77}.
To do so, we define a convex functional $\varphi : \boldsymbol{V}_{0}^{*}\to [0, +\infty ]$ by
\begin{equation}
	\label{varp}
	\varphi (\boldsymbol{z}):=\begin{cases}
	\displaystyle \int_{\Omega }\widehat{\beta }(z)dx
	+\int_{\Gamma }\widehat{\beta }(z_{\Gamma })d\Gamma 
	& {\rm if} 
	\quad \boldsymbol{z}\in \boldsymbol{H}_{0}, 
	\widehat{\beta }(z)\in L^{1}(\Omega ), \widehat{\beta }(z_{\Gamma })\in L^{1}(\Gamma ), \\
	+\infty & {\rm otherwise}. \\
\end{cases} 
\end{equation}
Note that the assumption of the growth condition {\rm (A2)} plays an important role for the 
lower semicontinuity on $\boldsymbol{V}_{0}^{*}$ of $\varphi $.

\paragraph{Lemma 3.1.}
{\it The proper convex functional $\varphi:\boldsymbol{V}_0^* \to [0,+\infty ]$ is lower semicontinuous 
on $\boldsymbol{V}_{0}^{*}$.
}

\paragraph{Proof}
It is enough to show that 
the level set $[\varphi \le \lambda ]:=\{ \boldsymbol{z} \in \boldsymbol{V}_0^* : \varphi (\boldsymbol{z}) \le \lambda \}$ is closed 
in $\boldsymbol{V}_0^*$ for all $\lambda \in \mathbb{R}$ (see, e.g., \cite[p.70, Proposition 2.5]{BP12}).
We first take any $\{ \boldsymbol{z}_n \} _{n \in \mathbb{N}} \subset [\varphi \le \lambda ]$ with 
$\boldsymbol{z}_n \to \boldsymbol{z}$ in $\boldsymbol{H}_0$ as 
$n\to +\infty $. 
$\widehat{\beta }$ is now lower semicontinuous on $\mathbb{R}$. Therefore,
by applying the {F}atou lemma to subsequences $\{ z_{n_k} \}_{k \in \mathbb{N}}$ 
and $\{ z_{\Gamma, n_k} \}_{k \in \mathbb{N}}$, 
which respectively converge to $z$ and $z_\Gamma $ almost everywhere, 
we see that 
$\varphi (\boldsymbol{z}) \le \liminf_{k \to \infty } \varphi (\boldsymbol{z}_{n_k}) \le \lambda$; 
i.e., $[\varphi \le \lambda ]$ is closed with respect to the topology of $\boldsymbol{H}_0$. 
Second, from the convexity of $\varphi$, we see that 
$[\varphi \le \lambda ]$ is closed with respect to the weak topology of $\boldsymbol{H}_0$ 
(see, e.g., \cite[p.72, Proposition 2.10]{BP12}).
We finally take any $\{ \boldsymbol{z}_n \} _{n \in \mathbb{N}} \subset [\varphi \le \lambda ]$ with 
$\boldsymbol{z}_n \to \boldsymbol{z}$ in $\boldsymbol{V}_0^*$ as 
$n\to +\infty $. 
In this case, from the assumption of growth condition {\rm (A2)}, 
we can take a bounded subsequence 
$\{ \boldsymbol{z}_{n_k} \}_{k \in \mathbb{N}}$ in $\boldsymbol{H}_0$ 
such that $\boldsymbol{z}_{n_k} \to \boldsymbol{z}$ weakly in $\boldsymbol{H}_0$ as $k \to +\infty $. 
We thus conclude that $\boldsymbol{z} \in [\varphi \le \lambda ]$. \hfill $\Box$\\

We now define the projection $\boldsymbol{P}:\boldsymbol{H} \to \boldsymbol{H}_0$ by 
\begin{equation}
	\boldsymbol{P} \boldsymbol{z}:=\boldsymbol{z} - m(\boldsymbol{z}) \boldsymbol{1} \quad 
	{\rm for~all~}\boldsymbol{z} \in \boldsymbol{H},
	\nonumber \label{P}
\end{equation}
where $\boldsymbol{1}:=(1,1)$.

\subsection{Approximate problem for (P)}

We next consider an approximate problem to show the existence of a weak solution to (P). 
For each $\lambda >0$, we define the Moreau--Yosida regularization $\widehat{\beta }_{\lambda }$ of $\widehat{\beta }: \mathbb{R}\to \mathbb{R}$ by
\begin{equation}
	\widehat{\beta }_{\lambda }(r)
	:=\inf _{s\in \mathbb{R}}
	\left\{\frac{1}{2\lambda }|r-s|^{2}+\widehat{\beta }(s)\right\}
	=\frac{1}{2\lambda }
	\bigl| r-J_{\lambda }(r) \bigr|^{2}
	+\widehat{\beta }\bigl( J_{\lambda }(r) \bigr)
	\label{MY}
\end{equation}
for all $r\in \mathbb{R}$, 
where the resolvent operator 
$J_{\lambda }: \mathbb{R}\to \mathbb{R}$ of $\beta $ 
is given by $J_{\lambda }(r):=(I+\lambda \beta )^{-1}r$. 
We also define 
\begin{equation*}
	\varphi _{\lambda }(\boldsymbol{z}):=
	\begin{cases}
	\displaystyle \int_{\Omega }\widehat{\beta }_{\lambda }(z)dx
	+\int_{\Gamma }\widehat{\beta }_{\lambda }(z_{\Gamma })d\Gamma \quad 
	& {\rm if~}\boldsymbol{z}\in \boldsymbol{H}_{0}, \\
	+\infty \quad & {\rm otherwise}. 
	\end{cases} 
\end{equation*}
Then, for each $\lambda >0$, $\varphi _{\lambda }$ is a proper lower semicontinuous convex function on $\boldsymbol{V}_{0}^{*}$. 
We now give the representation of subdifferential 
operator $\partial _{\boldsymbol{V}_{0}^{*}}\varphi _{\lambda }$ by the next lemma.

\paragraph{Lemma 3.2.}
{\it 
For any $\boldsymbol{z}\in \boldsymbol{H}_{0}$, the following equivalence holds:
$\boldsymbol{z}^* \in \partial _{\boldsymbol{V}_{0}^{*}}\varphi _{\lambda }(\boldsymbol{z})$ in 
$\boldsymbol{V}_0^*$ if and only if $\boldsymbol{\beta }_{\lambda }(\boldsymbol{z}):=(\beta _{\lambda }(z), \beta _{\lambda }(z_{\Gamma })) 
\in \boldsymbol{V}$ and 
\begin{equation*}
	\boldsymbol{z}^* =\boldsymbol{F}\boldsymbol{P}\boldsymbol{\beta }_{\lambda }(\boldsymbol{z}) 
	\quad {\it in~} \boldsymbol{V}_{0}^{*}. 
\end{equation*}
That is to say, $\partial _{\boldsymbol{V}_{0}^{*}}\varphi _{\lambda }(\boldsymbol{z})$ is a singleton.
}

\paragraph{Proof}
For each fixed $\boldsymbol{z}\in D(\varphi _{\lambda })=\boldsymbol{H}_{0}$, 
we set $\boldsymbol{z}^{*} \in \partial _{\boldsymbol{V}_{0}^{*}}\varphi _{\lambda }(\boldsymbol{z})\in \boldsymbol{V}_{0}^{*}$. 
We see from the definition of the subdifferential that $$(\boldsymbol{z}^{*}, \tilde{\boldsymbol{z}}-\boldsymbol{z})_{\boldsymbol{V}_{0}^{*}}\leq \varphi _{\lambda }(\tilde{\boldsymbol{z}})-\varphi _{\lambda }(\boldsymbol{z}) \quad {\rm for~all~}\tilde{\boldsymbol{z}}\in \boldsymbol{V}_{0}^{*}. $$
Now, for each $\delta \in (0, 1]$ and $\bar{\boldsymbol{z}}\in \boldsymbol{H}_{0}$, 
taking $\tilde{\boldsymbol{z}}:=\boldsymbol{z}+\delta \bar{\boldsymbol{z}}
\in \boldsymbol{H}_{0}$ in the above, we have
\begin{equation}
	\label{lem2-1}
	(\boldsymbol{z}^{*}, \delta \bar{\boldsymbol{z}})_{\boldsymbol{V}_{0}^{*}}\leq \int_{\Omega }\widehat{\beta }_{\lambda }(z+\delta \bar{z})dx-\int_{\Omega }\widehat{\beta }_{\lambda }(z)dx+\int_{\Gamma }\widehat{\beta }_{\lambda }(z_{\Gamma }+\delta \bar{z}_{\Gamma })d\Gamma -\int_{\Gamma }\widehat{\beta }_{\lambda }(z_{\Gamma })d\Gamma . 
\end{equation}
Now, according to the intermediate value theorem, 
there exist $\xi : \Omega \to \mathbb{R}$ 
between $z$ and $z+\delta \bar{z}$ a.e.\ in $\Omega $ and 
$\xi _{\Gamma }: \Gamma \to \mathbb{R}$ between $z_{\Gamma }$ and $z_{\Gamma }+\delta \bar{z}_{\Gamma }$ a.e.\ on $\Gamma$ such that
\begin{gather*}
	\frac{\widehat{\beta }_{\lambda }(z+\delta \bar{z})
	-\widehat{\beta }_{\lambda }(z)}{\delta }
	=\beta_{\lambda }(\xi )\bar{z} \quad {\rm a.e.~in}~\Omega, 
	\\
	\frac{\widehat{\beta }_{\lambda }(z_{\Gamma }+\delta \bar{z}_{\Gamma })
	-\widehat{\beta }_{\lambda }(z_{\Gamma })}{\delta }
	=\beta_{\lambda }(\xi _{\Gamma })\bar{z}_{\Gamma } \quad {\rm a.e.~on}~\Gamma.
\end{gather*}
We thus deduce that 
\begin{align*}
	\left|
	\frac{\widehat{\beta }_{\lambda }(z+\delta \bar{z})-\widehat{\beta }_{\lambda }(z)}{\delta }
	\right|
	&=|\beta _{\lambda }(\xi )-\beta _{\lambda }(0)||\bar{z}| \\
	&\leq \frac{1}{\lambda }|\xi -0||\bar{z}| \\
	&\leq \frac{1}{\lambda }(|z|+\delta |\bar{z}|)|\bar{z}|
\end{align*}
a.e.\ in $\Omega $, where the {L}ipschitz continuity of $\beta _{\lambda }$ with the {L}ipschitz constant 
$1/\lambda $ is used. 
Now, letting $\delta $ tend to zero, we obtain $\xi \to z$ a.e.\ in 
$\Omega $, $\beta _{\lambda }(\xi )\to \beta _{\lambda }(z)$ a.e.\ in $\Omega $. 
From the {L}ebesgue dominated convergence theorem, it follows that
\begin{equation*}
	\lim _{\delta \to 0}\int_{\Omega }\frac{\widehat{\beta }_{\lambda }(z+\delta \bar{z})
	-\widehat{\beta }_{\lambda }(z)}{\delta }dx
	=\int_{\Omega }\beta _{\lambda }(z)\bar{z}dx= \bigl( \beta _{\lambda }(z), \bar{z} \bigr)_{H}. 
\end{equation*}
Similarly, 
\begin{equation*}
	\lim _{n\to \infty }\int_{\Gamma }\frac{\widehat{\beta }_{\lambda }
	(z_{\Gamma }+\delta \bar{z}_{\Gamma })-\widehat{\beta }_{\lambda }(z_{\Gamma })}{\delta }d\Gamma
	 =\bigl( \beta _{\lambda }(z_{\Gamma }), \bar{z}_{\Gamma } \bigr)_{H_{\Gamma }}. 
\end{equation*}
Thus, 
through dividing by $\delta >0$ in \eqref{lem2-1} and letting $\delta $ tend to zero, we infer that 
\begin{align*}
	(\boldsymbol{z}^{*}, \bar{\boldsymbol{z}})_{\boldsymbol{V}_{0}^{*}}
	& \leq 
	\bigl(
	\beta _{\lambda }(z), \bar{z} \bigr)_{H} + 
	\bigl( 
	\beta _{\lambda }(z_{\Gamma }), \bar{z}_{\Gamma } 
	\bigr)_{H_{\Gamma }} \\
	&= \bigl( 
	\boldsymbol{\beta }_{\lambda }(\boldsymbol{z}), \bar{\boldsymbol{z}} 
	\bigr)_{\boldsymbol{H}} \\
	&= \bigl( 
	\boldsymbol{P}\boldsymbol{\beta }_{\lambda }(\boldsymbol{z}), \bar{\boldsymbol{z}}
	\bigr)_{\boldsymbol{H}_{0}} \quad {\rm for~all~}\bar{\boldsymbol{z}}\in \boldsymbol{H}_{0}. 
\end{align*}
Next, taking $\tilde{\boldsymbol{z}}:=\boldsymbol{z}-\delta \bar{\boldsymbol{z}}$, we see that
$(\boldsymbol{z}^{*}, \bar{\boldsymbol{z}})_{\boldsymbol{V}_{0}^{*}}
\geq  (\boldsymbol{P}\boldsymbol{\beta }_{\lambda }(\boldsymbol{z}), \bar{\boldsymbol{z}})_{\boldsymbol{H}_{0}}$ 
for all $\bar{\boldsymbol{z}}\in \boldsymbol{H}_{0}$. 
That is to say, we have $(\boldsymbol{z}^{*}, \bar{\boldsymbol{z}})_{\boldsymbol{V}_{0}^{*}}=(\boldsymbol{P}\boldsymbol{\beta }_{\lambda }(\boldsymbol{z}),\bar{\boldsymbol{z}})_{\boldsymbol{H}_{0}} \ {\rm for~all~}\bar{\boldsymbol{z}}\in \boldsymbol{H}_{0}$. 
This implies that 
$\boldsymbol{F}^{-1}\boldsymbol{z}^{*}
=\boldsymbol{P}\boldsymbol{\beta }_{\lambda }(\boldsymbol{z})$ in $\boldsymbol{H}_{0}$, 
that is, in 
$\boldsymbol{V}_{0}$ by comparison. 
We therefore get $\boldsymbol{\beta }_{\lambda }(\boldsymbol{z})
\in \boldsymbol{V}$ and 
$\boldsymbol{z}^{*}=\boldsymbol{F}\boldsymbol{P}\boldsymbol{\beta }_{\lambda }(\boldsymbol{z})$ in $\boldsymbol{V}_{0}^{*}$. 
Meanwhile, 
if $\boldsymbol{\beta }_{\lambda }(\boldsymbol{z})
\in \boldsymbol{V}$, then $\boldsymbol{P} \boldsymbol{\beta }_{\lambda }(\boldsymbol{z})
\in \boldsymbol{V}_0$ and 
\begin{align}
	\bigl( \boldsymbol{F}\boldsymbol{P}\boldsymbol{\beta }_{\lambda }(\boldsymbol{z}), 
	\tilde{\boldsymbol{z}} - \boldsymbol{z} \bigr)_{\boldsymbol{V}_0^*} 
	& = \bigl \langle  
	\tilde{\boldsymbol{z}} - \boldsymbol{z}, 
	\boldsymbol{P}\boldsymbol{\beta }_{\lambda }(\boldsymbol{z})
	\bigr \rangle _{\boldsymbol{V}_0^*, \boldsymbol{V}_0} \nonumber \\
	& = \bigl( \tilde{\boldsymbol{z}} - \boldsymbol{z}, 
	\boldsymbol{\beta }_{\lambda }(\boldsymbol{z}) \bigr)_{\boldsymbol{H}} \nonumber \\
	& \le \varphi_\lambda (\tilde{\boldsymbol{z}}) - \varphi _\lambda (\boldsymbol{z}) 
	\quad {\rm for~all~} \tilde{\boldsymbol{z}} \in \boldsymbol{H}_0,
	\label{sub}
\end{align}
because $ \tilde{\boldsymbol{z}} - \boldsymbol{z} \in \boldsymbol{H}_0$. 
If $\tilde{\boldsymbol{z}} \in \boldsymbol{V}_0^* \setminus \boldsymbol{H}_0$, then 
$ \varphi_\lambda (\tilde{\boldsymbol{z}}) =+\infty $, and 
\eqref{sub} thus holds for all $\tilde{\boldsymbol{z}} \in \boldsymbol{V}_0^*$. This gives us 
$\boldsymbol{z}^* \in \partial _{\boldsymbol{V}_{0}^{*}}\varphi _{\lambda }(\boldsymbol{z})$ in 
$\boldsymbol{V}_0^*$. 
\hfill $\Box$ \\

For each $\lambda \in (0,1]$, applying the abstract theory of {B}r\'ezis (see \cite{Bre73}), 
we see that for each 
$\boldsymbol{f}$ and 
$\boldsymbol{u}_{0}$ satisfying (A3) and (A4), 
there exists a unique function $\boldsymbol{u}_{\lambda }\in H^{1}(0, T;\boldsymbol{V}_{0}^{*})\cap L^{\infty }(0, T;D(\varphi _{\lambda }))$ such that
\begin{equation}
	\begin{cases}
	\boldsymbol{u}_{\lambda }'(t)+\partial _{\boldsymbol{V}_{0}^{*}} \varphi _{\lambda }
	\bigl( \boldsymbol{u}_{\lambda }(t) \bigr)
	=\boldsymbol{f}(t) \quad {\rm in}~\boldsymbol{V}_{0}^{*} \quad {\rm for~a.a.~} t\in (0, T), \\
	\boldsymbol{u}_{\lambda }(0)=\boldsymbol{u}_{0} \quad {\rm in}~\boldsymbol{V}_{0}^{*}. 
	\end{cases} 
	\label{eel}
\end{equation}
From Lemma~3.2, it follows that 
$\boldsymbol{\beta }_{\lambda }(\boldsymbol{u}_\lambda (t))
\in \boldsymbol{V}$ for a.a.\ $t \in (0,T)$ and 
\begin{align*}
	\boldsymbol{f}(t)-\boldsymbol{u}_{\lambda }'(t) 
	& = \partial _{\boldsymbol{V}_{0}^{*}} \varphi _{\lambda }
	\bigl( \boldsymbol{u}_{\lambda }(t) \bigr) \nonumber \\
	& = 
	\boldsymbol{F}\boldsymbol{P}\boldsymbol{\beta }_{\lambda }
	\bigl( 
	\boldsymbol{u}_{\lambda }(t)
	\bigr) \quad {\rm in}~\boldsymbol{V}_{0}^{*}. 
	\label{eel2}
\end{align*}
This yields 
\begin{equation}
	\bigl\langle \boldsymbol{u}_{\lambda }'(t), \boldsymbol{z}
	\bigr\rangle _{\boldsymbol{V}_{0}^{*}, \boldsymbol{V}_{0}}
	+a\bigl( 
	\boldsymbol{\beta }_{\lambda }
	\bigl( \boldsymbol{u}_{\lambda }(t) \bigr), 
	\boldsymbol{z} \bigr) 
	= \bigl( \boldsymbol{f}(t), \boldsymbol{z} \bigr)_{\boldsymbol{H}_{0}} 
	\quad {\rm for~all~}\boldsymbol{z}\in \boldsymbol{V}_{0},
	\label{eel3}
\end{equation}
for a.a.\ $t \in (0,T)$.

\subsection{Uniform estimates}

In this subsection, we obtain the uniform estimates, independent of $\lambda $, 
to prove the suitable convergence.

\paragraph{Lemma 3.3.}
{\it 
There exists a positive constant $M_{1}$, independent of $\lambda \in (0, 1]$, such that
\begin{equation*}
	\frac{1}{2}\int_{0}^{t} 
	\bigl| 
	\boldsymbol{u}_{\lambda }'(s)
	\bigr|_{\boldsymbol{V}_{0}^{*}}^{2}ds 
	+ \int_{\Omega }\widehat{\beta }_{\lambda }
	\bigl( u_{\lambda }(t) \bigr)dx
	+ \int_{\Gamma }\widehat{\beta }_{\lambda }
	\bigl( u_{\Gamma , \lambda }(t) \bigr) d\Gamma \leq M_{1}
\end{equation*}
for all $t\in [0, T]$. }

\paragraph{Proof}
For a.a.\ $s\in (0, T)$, 
we have that 
\begin{equation*}
	\bigl(\boldsymbol{u}_{\lambda }'(s), \partial _{\boldsymbol{V}_{0}^{*}}
	\varphi _{\lambda }\bigl( 
	\boldsymbol{u}_{\lambda }(s) \bigr) \bigr)_{\boldsymbol{V}_{0}^{*}}
	=\frac{d}{ds}\varphi _{\lambda } \bigl( 
	\boldsymbol{u}_{\lambda }(s) \bigr). 
\end{equation*}
Hence, we deduce from \eqref{eel} that
\begin{align*}
	\bigl|\boldsymbol{u}_{\lambda }'(s) \bigr|_{\boldsymbol{V}_{0}^{*}}^{2} 
	& = \bigl( \boldsymbol{u}_{\lambda }'(s), \boldsymbol{f}(s)
	- \partial _{\boldsymbol{V}_{0}^{*}}\varphi _{\lambda }
	\bigl( 
	\boldsymbol{u}_{\lambda }(s) \bigr) \bigr)_{\boldsymbol{V}_{0}^{*}} \\
	& = \bigl\langle  \boldsymbol{u}_{\lambda }'(s), \boldsymbol{F}^{-1}
	\boldsymbol{f}(s)\bigr\rangle _{\boldsymbol{V}_{0}^{*}, \boldsymbol{V}_{0}}
	-\frac{d}{ds}\varphi _{\lambda } \bigl( \boldsymbol{u}_{\lambda }(s) \bigr) . 
\end{align*}
Now, integrating over $(0, t)$ with respect to $s$, we obtain 
\begin{equation*}
	\int_{0}^{t}\bigl|\boldsymbol{u}_{\lambda }'(s) \bigr|_{\boldsymbol{V}_{0}^{*}}^{2}ds
	+\varphi _{\lambda } \bigl( \boldsymbol{u}_{\lambda }(t) \bigr)
	\leq \varphi (\boldsymbol{u}_{0})+\int_{0}^{t}\bigl|\boldsymbol{u}_{\lambda }'(s)\bigr|_{\boldsymbol{V}_{0}^{*}}\bigl|\boldsymbol{F}^{-1}\boldsymbol{f}(s)\bigr|_{\boldsymbol{V}_{0}}ds 
\end{equation*}
for all $t \in [0,T]$. 
Then, using the {Y}oung inequality and taking 
\begin{equation*}
	M_{1}:=\bigl| \widehat{\beta }(u_{0}) \bigr|_{L^{1}(\Omega )}
	+\bigl| \widehat{\beta }(u_{0\Gamma }) \bigr|_{L^{1}(\Gamma )}
	+\frac{1}{2}|\boldsymbol{f}|_{L^{2}(0, T;\boldsymbol{H}_{0})}^2,
\end{equation*} 
we get the conclusion. 
\hfill $\Box $

\paragraph{Lemma 3.4.}
{\it 
There exist 
a value $\bar{\lambda }\in (0, 1]$ and 
a positive constant $M_{2}$ independent of $\lambda \in (0, \bar{\lambda }]$, 
such that
\begin{equation}
	\bigl| \boldsymbol{u}_{\lambda }(t) \bigr|_{\boldsymbol{H}_{0}}^{2} \leq M_{2}
	\label{lem3.4}
\end{equation}
for all $t\in [0, T]$ and $\lambda \in (0, \bar{\lambda }]$. }

\paragraph{Proof}
By virtue of \eqref{MY} with (A2), we have that
\begin{equation*}
	\widehat{\beta }_{\lambda } 
	\bigl( 
	u_{\lambda }(t) \bigr) 
	\geq \frac{1}{2\lambda }
	\bigl| 
	u_{\lambda }(t)-J_{\lambda } \bigl( u_{\lambda }(t) \bigr)
	\bigr|^{2}
	+c_{1}
	\bigl| 
	J_{\lambda }
	\bigl(u_{\lambda }(t)
	\bigr)
	\bigr|^{2}-c_{2},
\end{equation*}
for all $t \in [0,T]$. 
We now set $\bar{\lambda }:=\min \left\{1, 1/(2c_{1})\right\}$. Then, 
for each $\lambda \in (0, \bar{\lambda }]$, we have
$\lambda \leq \bar{\lambda }\leq 1/(2c_{1})$; i.e., 
$1/(2\lambda )\geq c_{1}$. It follows from Lemma~3.3 that
\begin{align*}
	M_{1} & \geq \int_{\Omega } \widehat{\beta }_{\lambda } \bigl( u_{\lambda }(t) \bigr)dx \\
	& \geq c_{1}\int_{\Omega } \left\{ 
	\bigl| u_{\lambda }(t)-J_{\lambda } \bigl (u_{\lambda }(t) \bigr) \bigr |^{2} 
	+ \bigl| J_{\lambda } \bigl (u_{\lambda }(t) \bigr) \bigr |^{2} \right\} dx 
	- c_{2}|\Omega | \\
	& \geq \frac{c_1}{2} \int_{\Omega } 
	\bigl| u_{\lambda }(t) \bigr|^{2}dx-c_{2}|\Omega |. 
\end{align*}
This yields 
\begin{equation*}
	\bigl| 
	u_{\lambda }(t) \bigr|_{H}^{2} \leq 
	\frac{2}{c_{1}} 
	\bigl( M_{1}+c_{2}|\Omega | \bigr) 
	\quad {\rm for~all~}t\in [0, T]. 
\end{equation*}
Similarly, 
\begin{equation*}
	\bigl| 
	u_{\Gamma , \lambda }(t)
	\bigr|_{H_{\Gamma }}^{2}\leq 
	\frac{2}{c_{1}} 
	\bigl( M_{1}+c_{2}|\Gamma | \bigr) 
	\quad {\rm for~all~}t\in [0, T]. 
\end{equation*}
Thus, setting $M_{2}:=(2/c_{1})(2M_{1}+c_{2}(|\Omega |+|\Gamma |))$, we obtain \eqref{lem3.4}. 
\hfill 
$\Box $

\paragraph{Lemma 3.5.}
{\it 
There exist positive constants $M_{3}$ and $M_{4}$, independent of $\lambda \in (0, 1]$, such that
\begin{equation*}
	\int_{0}^{t} \bigl| 
	\boldsymbol{P}\boldsymbol{\beta }_{\lambda }
	\bigl(\boldsymbol{u}_{\lambda }(s) \bigr) \bigr |_{\boldsymbol{V}_{0}}^{2}ds\leq M_{3},
	\quad 
	\bigl| 
	m\bigl (\boldsymbol{\beta }_{\lambda }\bigl( \boldsymbol{u}_{\lambda }(t) \bigr)
	\bigr ) \bigr |\leq M_{4}
\end{equation*}
for all $t\in [0, T]$. }

\paragraph{Proof}
For all $t\in [0, T]$, from \eqref{eel} with Lemma~3.3, we deduce that
\begin{align*}
	\int_{0}^{t} 
	\bigl| 
	\partial _{\boldsymbol{V}_{0}^{*}}
	\varphi _{\lambda } \bigl(\boldsymbol{u}_{\lambda }(s)
	\bigr) \bigr |_{\boldsymbol{V}_{0}^{*}}^{2}ds
	& =\int_{0}^{t}
	\bigl| 
	\boldsymbol{f}(s)-\boldsymbol{u}_{\lambda }'(s) 
	\bigr|_{\boldsymbol{V}_{0}^{*}}^{2}ds \\ 
	&\leq 2 \int_{0}^{T} \bigl| 
	\boldsymbol{f}(s) \bigr|_{\boldsymbol{V}_{0}^{*}}^{2}ds
	+2\int_{0}^{T}
	\bigl| 
	\boldsymbol{u}_{\lambda }'(s) 
	\bigr|_{\boldsymbol{V}_{0}^{*}}^{2}ds \\
	&\leq 2|\boldsymbol{f}|_{L^{2}(0, T;\boldsymbol{H}_{0})}^{2}+4M_{1}. 
\end{align*}
Now, by setting $M_{3}:=2|\boldsymbol{f}|_{L^{2}(0, T;\boldsymbol{H}_{0})}^{2}+4M_{1}$ we infer from Lemma~3.2 that
\begin{align*}
	\int_{0}^{t}
	\bigl| 
	\boldsymbol{P}\boldsymbol{\beta }_{\lambda } 
	\bigl( \boldsymbol{u}_{\lambda }(s) \bigr)
	\bigr|_{\boldsymbol{V}_{0}}^{2}ds
	& =
	\int_{0}^{t} 
	\bigl| 
	\boldsymbol{F}\boldsymbol{P}\boldsymbol{\beta }_{\lambda } 
	\bigl( \boldsymbol{u}_{\lambda } (s) \bigr) 
	\bigr|_{\boldsymbol{V}_{0}^{*}}^{2}ds \\
	& = \int_{0}^{t} 
	\bigl| 
	\partial _{\boldsymbol{V}_{0}^{*}}\varphi _{\lambda } \bigl( 
	\boldsymbol{u}_{\lambda }(s)
	\bigr) 
	\bigr|_{\boldsymbol{V}_{0}^{*}}^{2}ds \\
	& \leq M_{3}
\end{align*}
for all $t \in [0,T]$.
To prove the second estimate, we set 
\begin{gather*}
	\Omega _{1}=\Omega _{1}(u_{\lambda })
	:=\bigl\{ x\in \Omega ; |u_{\lambda }(x)|\leq M_{0} \bigr\}, 
	\quad \Omega _{2}=\Omega _{2}(u_{\lambda })
	:=\bigl \{x\in \Omega ; |u_{\lambda }(x)|>M_{0} \bigr \}, \\
	\Gamma _{1}=\Gamma _{1}(u_{\Gamma , \lambda })
	:=\bigl \{x\in \Gamma ; |u_{\Gamma , \lambda }(x)|\leq M_{0} \bigr \}, 
	\quad \Gamma _{2}=\Gamma _{2}(u_{\Gamma , \lambda })
	:=\bigl \{x\in \Gamma ; |u_{\Gamma , \lambda }(x)|>M_{0} \bigr \},  
\end{gather*} 
and positive constant $c^*:=\max \{\beta (M_{0}), -\beta (-M_{0})\}$. 
For all $t\in [0, T]$, we infer that 
\begin{align*}
	& \bigl| 
	m\bigl (\boldsymbol{\beta }_{\lambda }\bigl( \boldsymbol{u}_{\lambda }(t) \bigr)
	\bigr ) \bigr | \\
	& \quad = \frac{1}{|\Omega |+|\Gamma |} 
	\left|\int_{\Omega } \beta _{\lambda }\bigl( 
	u_{\lambda }(t) \bigr)dx 
	+ \int_{\Gamma }\beta _{\lambda } \bigl( 
	u_{\Gamma , \lambda }(t) \bigr) d\Gamma \right| \\
	& \quad \leq \frac{1}{|\Omega |+|\Gamma |}
	\left\{ c^{*} \bigl(|\Omega | + |\Gamma | \bigr)
	+ c_0 \left|\int_{\Omega _{2}}u_{\lambda }(t)dx 
	+ \int_{\Gamma _{2}} u_{\Gamma, \lambda }(t) d \Gamma \right|  
	+ c_0' \bigl(|\Omega | + |\Gamma | \bigr) \right\} . 
\end{align*}
Now, from total mass conservation we have 
\begin{align*}
	\left| \int_{\Omega _{2}} u_{\lambda }(t)dx
	+ \int_{\Gamma _{2}}  u_{\Gamma , \lambda }(t)d\Gamma
	\right| 
	& = 
	\left|  
	- \left( \int_{\Omega _{1}}  u_{\lambda }(t)dx
	+ \int_{\Gamma _{1}}  u_{\Gamma , \lambda }(t)d\Gamma 
	\right) 
	\right| \\
	& \le  M_{0}(|\Omega |+|\Gamma |). 
\end{align*}
We thus have 
\begin{equation*}
	\bigl| 
	m\bigl (\boldsymbol{\beta }_{\lambda }\bigl( \boldsymbol{u}_{\lambda }(t) \bigr)
	\bigr )
	\bigr| 
	\leq c^{*}+c_{0}M_0 + c_0'=:M_{4}
\end{equation*}
for all $t \in [0,T]$. 
\hfill  $\Box $ \\

\paragraph{Lemma 3.6.}
{\it 
There exists a positive constant $M_{5}$, independent of $\lambda \in (0,1]$, such that
\begin{align*}
	\int_{0}^{t}
	\bigl| \boldsymbol{\beta }_{\lambda }
	\bigl( 
	\boldsymbol{u}_{\lambda }(s) 
	\bigr) \bigr|_{\boldsymbol{V}}^{2}ds\leq M_{5}
\end{align*}
for all $t\in [0, T]$. }

\paragraph{Proof}
We consider that $a(\boldsymbol{z}, \boldsymbol{z})=|\boldsymbol{P}\boldsymbol{z}|_{\boldsymbol{V}_0}^2$
for all $\boldsymbol{z} \in \boldsymbol{V}$. 
From \eqref{PW0} with Lemma~3.5, we have
\begin{align*}
	\int_{0}^{t}
	\bigl| 
	\boldsymbol{\beta }_{\lambda }
	\bigl(\boldsymbol{u}_{\lambda }(s) \bigr) 
	\bigr|_{\boldsymbol{V}}^{2}ds
	& 
	\leq c_{{\rm P}}
	\left\{
	\int_{0}^{t} 
	\bigl| 
	\boldsymbol{P} \boldsymbol{\beta }_\lambda 
	\bigr( \boldsymbol{u}_\lambda (s) \bigr)
	\bigr|^2_{\boldsymbol{V}_0} ds 
	+
	\int_{0}^{t}
	\bigl| 
	m \bigl( 
	\boldsymbol{\beta }_{\lambda }(\boldsymbol{u}_{\lambda }(s) \bigr)
	\bigr|^2 ds \right\} \\
	&\leq c_{{\rm P}}(M_{3}+M_{4}^2T)=:M_{5}
\end{align*}
for all $t \in [0,T]$
\hfill $\Box $ \\

\subsection{Passage to the limit as $\lambda \to 0$}

In this subsection, we obtain the weak solution of (P) from the passage to the 
limit for the approximate problem.

\paragraph{Proof of Theorem 2.1.}
On the basis of the previous estimates in Lemmas~3.3, 3.4, and 3.6,
there exist a subsequence $\{\lambda _{k}\}_{k\in \mathbb{N}}$ with 
$\lambda _{k} \to 0$ as $k \to  \infty $ and limit functions 
$\boldsymbol{u}\in H^{1}(0, T; \boldsymbol{V}_{0}^{*})\cap L^{\infty }(0, T;\boldsymbol{H}_{0})$ 
and $\boldsymbol{\xi }\in L^{2}(0, T; \boldsymbol{V})$ such that 
\begin{gather}
	\boldsymbol{u}_{\lambda _{k}}\to \boldsymbol{u} 
	\quad {\rm weakly}~{\rm star~in}~H^{1}(0, T; \boldsymbol{V}_{0}^{*})\cap L^{\infty }(0, T; \boldsymbol{H}_{0}), \label{p-1} \\
	\boldsymbol{\beta }_{\lambda _{k}}(\boldsymbol{u}_{\lambda _{k}})\to \boldsymbol{\xi } \quad {\rm weakly~in}~L^{2}(0, T; \boldsymbol{V}) \label{p-3}
\end{gather}
as $k\to +\infty $. 
Now, from \eqref{p-1} and the {A}scoli--{A}rzela theorem (see, e.g., \cite{Sim87}), we see that 
there exists a subsequence (not relabeled) such that
\begin{equation}
	\label{p-4}
	\boldsymbol{u}_{\lambda _{k}}\to \boldsymbol{u} 
	\quad {\rm in}~C \bigl( [0, T]; \boldsymbol{V}_{0}^{*} \bigr)
\end{equation}
as $k\to +\infty $; i.e., 
$\boldsymbol{u}(0)=\boldsymbol{u}_{0}$ in $\boldsymbol{V}_{0}^{*}$. 
Now, from \eqref{eel3} and by letting $k \to \infty $, we obtain
\begin{equation}
	\bigl\langle  \boldsymbol{u}'(t), \boldsymbol{z}
	\bigr\rangle _{\boldsymbol{V}_{0}^{*}, \boldsymbol{V}_{0}}
	+a \bigl( \boldsymbol{\xi }(t), \boldsymbol{z} \bigr)
	= \bigl( \boldsymbol{f}(t), \boldsymbol{z} \bigr)_{\boldsymbol{H}_{0}}
	 \quad {\rm for~all~}\boldsymbol{z}\in \boldsymbol{V}_{0},
	 \label{lastw}
\end{equation}
for a.a.\ $t\in (0,T)$. 
To prove the main theorem, we show that 
$\xi (t)\in \beta(u(t))$ a.e.\ in 
$\Omega $ and 
$\xi_\Gamma  (t)\in \beta(u_\Gamma (t))$ a.e.\ on
$\Gamma $.
We now define two operators 
$\boldsymbol{B}, \boldsymbol{B}_{\lambda _{k}}: 
L^2(0,T;\boldsymbol{H}) \to L^2(0,T;\boldsymbol{H})$ by
\begin{gather*}
	\boldsymbol{B\eta }
	:= \bigl\{ 
	\boldsymbol{\xi}:=(\xi, \xi_{\Gamma })\in L^2(0,T;\boldsymbol{H}) \ : \
	\xi \in \beta (\eta ) \ {\rm a.e.\ in~}Q, \ \xi _{\Gamma }\in \beta (\eta_{\Gamma }) 
	\ {\rm a.e.\ on~}\Sigma \bigr\}, \\
	\boldsymbol{B}_{\lambda _{k}} \boldsymbol{\eta }
	:= \bigl( \beta _{\lambda _{k}}(\eta ), \beta _{\lambda _{k}}(\eta_\Gamma ) \bigr)
	\quad {\rm for~all~} \boldsymbol{\eta} \in L^2(0,T;\boldsymbol{H}).
\end{gather*}
Then, from the 
maximal monotonicity of $\beta $, we see that 
$\boldsymbol{B}$ and $\boldsymbol{B}_{\lambda _{k}}$ are 
maximal monotone operators on $L^2(0,T;\boldsymbol{H})$. 
Now, from \eqref{p-1} and \eqref{p-3} we already have 
\begin{gather*}
	\boldsymbol{u}_{\lambda _{k}}\to \boldsymbol{u} 
	\quad {\rm weakly~star~in~}L^{\infty }(0, T; \boldsymbol{H}), \\
	\boldsymbol{B}_{\lambda _{k}}\boldsymbol{u}_{\lambda _{k}}
	=\boldsymbol{\beta }_{\lambda _{k}}(\boldsymbol{u}_{\lambda _{k}})
	\to \boldsymbol{\xi } \quad {\rm weakly~in~}L^{2}(0, T; \boldsymbol{V})
\end{gather*}
as $k \to + \infty $. 
Moreover, we deduce from \eqref{p-4} that 
\begin{align*}
	\int_{0}^{T}
	\bigl( 
	\boldsymbol{B}_{\lambda _{k}} 
	\boldsymbol{u}_{\lambda _{k}}(s), \boldsymbol{u}_{\lambda _{k}}(s)
	\bigr )_{\boldsymbol{H}}ds 
	&= 
	\int_{0}^{T} \bigl\langle 
	\boldsymbol{u}_{\lambda _{k}}(s), 
	\boldsymbol{\beta }_{\lambda _{k}}
	\bigl( \boldsymbol{u}_{\lambda _{k}}(s) \bigr)
	\bigr\rangle _{\boldsymbol{V}^{*}, \boldsymbol{V}}ds \\
	&\to \int_{0}^{T} \bigl \langle 
	\boldsymbol{u}(s), \boldsymbol{\xi }(s) 
	\bigr \rangle _{\boldsymbol{V}^{*}, \boldsymbol{V}}ds \\
	&=\int_{0}^{T} \bigl (\boldsymbol{\xi }(s), \boldsymbol{u}(s) \bigr )_{\boldsymbol{H}}ds
\end{align*}
as $ k \to + \infty $. Therefore, 
by applying \cite[p.260, Lemma 7.1]{Ken07}, 
we deduce that $\boldsymbol{\xi } \in \boldsymbol{B}\boldsymbol{u}$ in 
$L^2(0,T;\boldsymbol{H})$. We finally check \eqref{uug} and 
\eqref{vi}. 
Firstly, 
it follows from \cite[Remark 2]{CF15} 
that the function $\boldsymbol{u}' \in L^2(0,T;\boldsymbol{V}_0^*)$ can be extended to 
$L^2(0,T;\boldsymbol{V}^*)$ by setting 
$\langle \boldsymbol{u}'(t), \boldsymbol{1} \rangle _{\boldsymbol{V}^*, \boldsymbol{V}}=0$; i.e.,
\begin{equation*}
	\bigl\langle \boldsymbol{u}'(t), \boldsymbol{z} \bigr\rangle_{\boldsymbol{V}^*, \boldsymbol{V}}
	:=  \bigl\langle \boldsymbol{u}'(t), \boldsymbol{P} \boldsymbol{z} \bigr\rangle_{\boldsymbol{V}^*, \boldsymbol{V}}
	\quad {\rm for~all~}\boldsymbol{z} \in \boldsymbol{V}.
\end{equation*}
We next see that $\boldsymbol{V}$ is a subspace of $V \times V_\Gamma$, 
and therefore, from the {H}ahn--{B}anach extension theorem, 
we can also extend $\boldsymbol{u}'(t)$ to $(V \times V_\Gamma)^*$; i.e.,  
$u \in H^1(0,T;V^*) \cap L^\infty (0,T;H)$ 
and $u_\Gamma \in H^1(0,T;V_\Gamma ^*) \cap L^\infty (0,T;H_\Gamma )$ 
with $\boldsymbol{u}=(u,u_\Gamma )\in H^1(0,T;\boldsymbol{V}^*) \cap L^\infty (0,T;\boldsymbol{H})$. 
Therefore, from \eqref{lastw} we obtain
\eqref{vi}. \hfill $\Box$

\section{Improvement}
\setcounter{equation}{0}

In this section, we consider an improvement to the main theorem.

\subsection{Improvement of the initial condition to a nonzero mean value}

The essential idea of the proof of Theorem~2.1 is the setting of 
suitable function spaces, or more precisely, the mean-zero function space $\boldsymbol{H}_0$. 
This idea comes from the treatment of the {C}ahn--{H}illiard system (see, e.g.,\ \cite{KNP95, CF15, Kub12}). 
Considering this idea, 
we can improve our assumption for the initial data to the general $\boldsymbol{H}$-function. 
In this subsection, 
we assume that
\begin{enumerate}
 \item[(A4)$^\prime $] $\boldsymbol{u}_{0}:=(u_{0}, u_{0\Gamma })
\in \boldsymbol{H}, \widehat{\beta }(u_{0})
\in L^{1}(\Omega )$ and $\widehat{\beta }(u_{0\Gamma })\in L^{1}(\Gamma )$. 
\end{enumerate}

We can then improve our main theorem as follows.
\paragraph{Theorem 4.1.}
{\it Under assumptions {\rm (A1)}--{\rm (A3)} and {\rm (A4)$^\prime $}, 
there exists a unique weak solution to the problem {\rm (P)}. }

\paragraph{Proof.} Let us consider a {C}auchy problem of the evolution equation: 
\begin{equation}
	\begin{cases}
	\boldsymbol{v}'(t)
	+\partial _{\boldsymbol{V}_{0}^{*}}
	\varphi_{m_0} \bigl( \boldsymbol{v}(t) \bigr) \ni \boldsymbol{f}(t) 
	\quad {\rm in}~\boldsymbol{V}_{0}^{*} \quad {\rm for~a.a.~}t\in (0, T), \\
	\boldsymbol{v}(0)=P\boldsymbol{u}_{0} \quad {\rm in}~\boldsymbol{H}_{0}, 
	\end{cases} 
	\label{ee2}
\end{equation}
where the convex functional $\varphi_{m_0} : \boldsymbol{V}_{0}^{*}\to [0, +\infty ]$ is defined by
\begin{equation}
	\label{varp2}
	\varphi_{m_0} (\boldsymbol{z}):=\begin{cases}
	\displaystyle \int_{\Omega }\widehat{\beta }(z+m_0)dx
	+ \int_{\Gamma }\widehat{\beta }(z_{\Gamma }+m_0)d\Gamma 
	\\
	\quad 
	{\rm if} 
	\quad \boldsymbol{z}\in \boldsymbol{H}_{0}, 
	\widehat{\beta }(z+m_0)\in L^{1}(\Omega ), \widehat{\beta }(z_{\Gamma }+m_0)\in L^{1}(\Gamma ), \\
	+\infty \quad {\rm otherwise}. \nonumber \\
\end{cases} 
\end{equation}
Then, this is also proper lower semicontinuous and convex on $\boldsymbol{V}_0^*$. 
Therefore, in the same way as for Theorem 2.1, we see that 
there exists a unique weak solution $(v,v_\Gamma, \xi, \xi_\Gamma)$ 
satisfying 
$\boldsymbol{v}\in H^{1}(0, T;\boldsymbol{V}_{0}^{*}) \cap L^{\infty }(0, T;\boldsymbol{H}_0)$
and $\boldsymbol{\xi } \in L^2(0,T;\boldsymbol{V})$ with 
$\xi \in \beta (v+m_0)$ a.e.\ in $Q$ and $\xi _\Gamma \in \beta (v_\Gamma +m_0)$ a.e.\ 
on $\Sigma $ such that \eqref{ee2} holds; i.e., $\boldsymbol{v}$ satisfies
\begin{gather*}
	\bigl\langle 
	v'(t), z
	\bigr\rangle _{V^{*}, V}+
	\bigl\langle  v'_{\Gamma }(t), z_{\Gamma }
	\bigr\rangle _{V_{\Gamma }^{*}, V_{\Gamma }}
	+\int_{\Omega }\nabla \xi (t)\cdot \nabla zdx
	+\int_{\Gamma }\nabla _{\Gamma }\xi _{\Gamma }(t)\cdot 
	\nabla _{\Gamma }z_{\Gamma }d\Gamma 
	\nonumber \\
	= 
	\int_{\Omega }
	f(t)zdx+\int_{\Gamma }f_{\Gamma }(t)z_{\Gamma }d\Gamma \quad 
	{\rm for~all~} \boldsymbol{z}\in \boldsymbol{V}, 
\end{gather*}
for a.a.\ $t\in (0, T)$. We now set $\boldsymbol{u}:=\boldsymbol{v}+m_0 \boldsymbol{1}$, and have 
$\boldsymbol{u}(0)=\boldsymbol{v}(0)+m_0 \boldsymbol{1}=\boldsymbol{u}_0$ in $\boldsymbol{H}$ and 
$\boldsymbol{u}'=\boldsymbol{v}'$. 
Thus, $(u,u_\Gamma, \xi, \xi _\Gamma)$ is our weak solution to problem (P). \hfill $\Box$

\subsection{Nonlinear diffusions of different $\beta $ and $\beta _\Gamma $}

In this subsection, we treat the different maximal monotone graphs 
$\beta $ and $\beta _{\Gamma }$ in $\Omega $ and on $\Gamma $, respectively. 
This implies that we can consider different nonlinearities of the diffusion in the bulk and on the boundary. 
We assume that

\begin{enumerate}
 \item[(A1)$^*$] $\beta, \beta _\Gamma  :\mathbb{R}\to 2^{\mathbb{R}}$ are maximal monotone graphs, which are
the subdifferentials $\beta =\partial _{\mathbb{R}}\widehat{\beta }$ and 
$\beta _\Gamma =\partial _\mathbb{R} \widehat{\beta }_\Gamma $
of some proper lower semicontinuous convex functions 
$\widehat{\beta }, \widehat{\beta }_\Gamma : \mathbb{R}\to [0, +\infty ]$ 
satisfying $\widehat{\beta }(0)=0$ and $\widehat{\beta }_\Gamma (0)=0$; 
 \item[(A2)$^*$] there exist positive constants $c_1$, $c_2$, $c_3$, and $c_4$ such that 
$\widehat{\beta }(r)\geq c_{1}r^{2}-c_{2}$ and 
$\widehat{\beta }_\Gamma (r)\geq c_{3}r^{2}-c_{4}$ for all $r\in \mathbb{R}$; 
 \item[(A4)$^*$] $\boldsymbol{u}_{0}:=(u_{0}, u_{0\Gamma })\in \boldsymbol{H}_0, 
\widehat{\beta }(u_{0})\in L^{1}(\Omega )$ and $\widehat{\beta }_\Gamma (u_{0\Gamma })\in L^{1}(\Gamma )$;
 \item[(A5)$^*$] there exist constants $c_0$, $M_0>0$, and $c_0', c_0'' \ge 0$
 such that
\begin{equation}
	\beta (r)
	=
	\begin{cases}
	c_0 r+c_0' & {\rm if~} r \ge  M_0, \\
	c_0 r-c_0' & {\rm if~} r \le  -M_0,
	\end{cases}
	\quad 
	\beta_\Gamma  (r)
	=
	\begin{cases}
	c_0 r+c_0'' & {\rm if~} r \ge  M_0, \\
	c_0 r-c_0'' & {\rm if~} r \le  -M_0.
	\end{cases}
	\nonumber 
\end{equation}

\end{enumerate}

We then obtain the same result.

\paragraph{Theorem 4.2.}
{\it Under assumptions {\rm (A1)$^*$}, {\rm (A2)$^*$}, {\rm (A3)}, {\rm (A4)$^*$}, and {\rm (A5)$^*$}
there exists a unique weak solution to the problem {\rm (P)} with 
\eqref{maximal} 
replaced by}
\begin{equation}
	\xi _{\Gamma }\in \beta_\Gamma (u_{\Gamma }), 
	\quad \xi _{|_{\Gamma }}=\xi _{\Gamma } 
	\quad {\it a.e.~on~}\Sigma. 
	\label{imp}
\end{equation}

\paragraph{Proof.}
We define a convex functional
\begin{gather}
	\label{chavar}
	\varphi (\boldsymbol{z}):=
	\begin{cases}
	\displaystyle \int_{\Omega }\widehat{\beta }(z)dx
	+\int_{\Gamma }\widehat{\beta }_{\Gamma }(z_{\Gamma })d\Gamma 
	& {\rm if} \quad \boldsymbol{z}\in \boldsymbol{H}_{0}, 
	\widehat{\beta }(z)\in L^{1}(\Omega ), 
	\widehat{\beta }_{\Gamma }(z_{\Gamma })\in L^{1}(\Gamma ), \\
	+\infty & {\rm otherwise}. \\
	\end{cases}
	\nonumber 
\end{gather}
(cf., \eqref{varp}). 
From growth conditions in {\rm (A2)$^*$}, 
we obtain coercivities of $\beta $ and $\beta _\Gamma $. 
We thus obtain 
surjectivities of $\beta $ and $\beta _\Gamma$; 
i.e., $R(\beta)=R(\beta _\Gamma )=\mathbb{R}$. 
Indeed, to obtain the trace condition in \eqref{imp}, we can use
the surjectivities of $\beta $ and $\beta _\Gamma $; 
otherwise, the trace condition makes no sense if $R(\beta ) \cap R(\beta _\Gamma )=\emptyset $. 
The entire proof of Theorem 4.2 is the same as that of Theorem 2.1, but with $\boldsymbol{\beta }$ replaced by $(\beta, \beta _\Gamma )$. 
\hfill $\Box$ \\

\paragraph{Remark 4.1.}
The characterization of 
the degenerate parabolic equation 
as the asymptotic limit of {C}ahn--{H}illiard systems has recently been discussed 
\cite{CF16, Fuk16, Fuk16b}. 
In the case of a dynamic boundary condition, for example, the existence result \cite{CF15} was used at the 
level of approximation. 
Comparing the proofs of Theorem~4.2 and \cite[Theorem~2.1]{Fuk16}, 
we readily see that {\rm (A5)} is too restrictive. 
Moreover, the growth condition in {(\rm A2)} has been improved \cite[Theorem~2.1]{Fuk16b} 
(see also the advantage of the Cahn--Hilliard approach \cite[Section~6]{CF16}). 
However, to treat different nonlinearities of $\beta $ and $\beta _\Gamma $, 
we will assume the condition (see \cite[p.419, {\rm (A6)}]{CF15}):
\begin{quote}
	There exist positive constants $\rho _1, \rho _2$ such that
\begin{equation}
	\bigl| \beta ^\circ (r) \bigr| \le \rho _1 \bigl| \beta _\Gamma ^\circ (r) \bigr|+\rho _2
	\quad {\rm for~all~} r\in \mathbb{R}.
	\label{cc}
\end{equation}
\end{quote}
Here, the minimal section $\beta ^\circ $ of $\beta $ is defined by 
$\beta ^\circ (r):=\{ r^* \in \beta (r) : |r^*|=\min _{s \in \beta (r)}|s| \}$ and 
the same definition applies to $\beta _\Gamma ^\circ $. 
Indeed, the dominated inequality \eqref{cc} is the same as 
\cite{CC13, CF15}, 
which gives us the same inequality at the level of the Yosida approximation \cite[p.19, Lemma 4.4]{CC13}. 
This dominated inequality provides 
suitable uniform estimates related to $\beta (u)$ and $\beta _\Gamma(u_\Gamma)$; 
more precisely, we can treat the estimate of $\beta _\Gamma (u_\Gamma )$ against $\beta (u_\Gamma )$ on the boundary  
(see \cite[Lemmas~4.3 and 4.4]{CF15}). 
However, if we apply the main theorem of the present paper, we do not need to assume such a dominated inequality because 
we do not treat directly 
the estimate of $\beta _\Gamma (u_\Gamma )$ against $\beta (u_\Gamma )$. 
This is one advantage of the abstract approach from the evolution equation  
to degenerate parabolic equations.

\section*{Acknowledgments}
We thank Edanz Group (www.edanzediting.com/ac) for editing a draft of this manuscript.


\begin{thebibliography}{99}


\bibitem{Aik93}
	T.\ {A}iki,
	\newblock Two-phase {S}tefan problems with dynamic boundary conditions, 
	\newblock Adv.\ Math.\ Sci.\ Appl., {\bf 2} (1993), 253--270. 

\bibitem{Aik95}
	T.\ {A}iki,
	\newblock Multi-dimensional {S}tefan problems with dynamic boundary conditions, 
	\newblock Appl.\ Anal., {\bf 56} (1995), 71--94. 

\bibitem{Aik96}
	T.\ {A}iki,
	\newblock Periodic stability of solutions to some degenerate parabolic equations with dynamic boundary conditions, 
	\newblock J.\ Math.\ Soc.\ Japan, {\bf 48} (1996), 37--59. 

\bibitem{AMTI06}
	F.\ {A}ndreu, J.\ M.\ {M}az{\'o}n, J.\ {T}oledo and N.\ {I}gbida,
	\newblock A degenerate elliptic-parabolic problem with nonlinear dynamical boundary conditions, 
	\newblock Interfaces Free Bound., {\bf 8} (2006), 447--479.

\bibitem{Bar10}
	V.\ {B}arbu,
	\newblock {\it Nonlinear differential equations of monotone types in {B}anach spaces}, 
	\newblock Springer, London, 2010. 

\bibitem{BP12} 
	\newblock V.\ {B}arbu and Th.\ {P}recupanu, 
	\newblock {\it Convexity and optimization in {B}anach space, 4th edition}, 
	\newblock Springer, Dordrecht, 2012. 

\bibitem{Bre73}
	H.\ {B}r\'ezis, 
	\newblock {\it Op\'erateurs maximaux monotones et semi-groupes de contractions dans les especes de {H}ilbert}, 
	\newblock North-Holland, Amsterdam, 1973.


\bibitem{CC13}
	L.\ {C}alatroni and P.\ {C}olli,
	\newblock Global solution to the {A}llen--{C}ahn equation with singular potentials and dynamic boundary conditions, 
	\newblock Nonlinear Anal., {\bf 79} (2013), 12--27.

\bibitem{CF15} 
	P.\ {C}olli and T.\ {F}ukao, 
	\newblock Equation and dynamic boundary condition of {C}ahn--{H}illiard type with singular potentials, 
	\newblock Nonlinear Anal., {\bf 127} (2015), 413--433.

\bibitem{CF16} 
	P.\ {C}olli and T.\ {F}ukao, 
	\newblock Nonlinear diffusion equations as asymptotic limits of {C}ahn--{H}illiard systems, 
	\newblock J.\ Differential Equations, {\bf 260} (2016), 6930--6959.


\bibitem{Dam77} 
	A.\ {D}amlamian, 
	\newblock Some results on the multi-phase {S}tefan problem, 
	\newblock Comm.\ Partial Differential Equations, {\bf 2} (1977), 1017--1044.

\bibitem{DK99}
	A.\ {D}amlamian and N.\ {K}enmochi,
	\newblock Evolution equations generated by subdifferentials in the dual space of $(H^1(\Omega))$, 
	\newblock Discrete Contin.\ Dyn.\ Syst., {\bf 5} (1999), 269--278.

\bibitem{Fuk16}
	T.\ {F}ukao, 
	\newblock Convergence of {C}ahn--{H}illiard systems to the {S}tefan problem with dynamic boundary conditions, 
	\newblock Asymptot.\ Anal., {\bf 99} (2016), 1--21. 

\bibitem{Fuk16b}
	T.\ {F}ukao, 
	\newblock {C}ahn--{H}illiard approach to some degenerate parabolic equations with dynamic boundary conditions, 
	\newblock in: L.\ Bociu, J-A. D\'esid\'eri and A.\ Habbal (Eds.), {\it System Modeling and Optimization}, Springer, Switzerland, 2016, pp.\ 282--291.


\bibitem{Gri09} 
	A.\ {G}rigor'yan,
	\newblock {\it Heat kernel and analysis on manifolds}, 
	\newblock American Mathematical Society, International Press, Boston, 2009.

\bibitem{Igb07}
	N.\ {I}gbida, 
	\newblock {H}ele-{S}haw type problems with dynamicalboundary conditions, 
	\newblock J.\ Math.\ Anal.\ Appl., {\bf 335} (2007), 1061--1078.

\bibitem{IK02}
	N.\ {I}gbida and M.\ {K}irane, 
	\newblock A degenerate diffusion problem with dynamical boundary conditions, 
	\newblock Math.\ Ann.\ {\bf 323} (2002), 377--396. 

\bibitem{Ken90}
	N.\ {K}enmochi, 
	\newblock Neumann problems for a class of nonlinear degenerate parabolic equations, 
	\newblock Differential Integral Equations, {\bf 2} (1990), 253--273.

\bibitem{Ken07}
	N.\ {K}enmochi,
	\newblock Monotonicity and compactness methods for nonlinear variational inequalities, 
	\newblock M.\ {C}hipot (Ed.), {\it Handbook of differential equations: 
	Stationary partial differential equations}, {\bf Vol.4}, 
	North-Holland, Amsterdam (2007), 203--298. 

\bibitem{KNP95}
	 N.\ {K}enmochi, M.\ {N}iezg\'odka and I.\ {P}aw\l ow, 
	\newblock Subdifferential operator approach to the {C}ahn--{H}illiard equation with constraint,
	\newblock J.\ Differential Equations, {\bf 117} (1995), 320--354.

\bibitem{Kub12} 
	M.\ {K}ubo,
	\newblock The {C}ahn--{H}illiard equation with time-dependent constraint, 
	\newblock Nonlinear Anal., {\bf 75} (2012), 5672--5685. 

\bibitem{KL05} 
	M.\ {K}ubo and Q.\ Lu,
	\newblock Nonlinear degenerate parabolic equations with {N}eumann boundary condition, 
	\newblock J.\ Math.\ Anal.\ Appl., {\bf 307} (2005), 232--244. 

\bibitem{KY17}
	S.\ {K}urima and T.\ {Y}okota,
	\newblock Monotonicity methods for nonlinear diffusion equations and their approximations with error estimates, 
	\newblock J.\ Differential Equations, {\bf 263} (2017), 2024--2050.



\bibitem{Sim87}
	 J.\ {S}imon, 
	\newblock Compact sets in the spaces $L^p(0,T;B)$, 
	\newblock Ann.\ Mat.\ Pura.\ Appl.\ (4), {\bf 146} (1987), 65--96.
\end{thebibliography}
\end{document}